\documentclass[12pt,reqno]{article}
\usepackage[usenames]{color}
\usepackage{amssymb}
\usepackage{graphicx}
\usepackage{amscd}
\usepackage{epsfig}
\usepackage[colorlinks=true,
linkcolor=webgreen,
filecolor=webbrown,
citecolor=webgreen]{hyperref}

\definecolor{webgreen}{rgb}{0,.5,0}
\definecolor{webbrown}{rgb}{.6,0,0}

\usepackage{color}
\usepackage{fullpage}
\usepackage{float}

\usepackage{psfig}
\usepackage{graphics,amsmath,amssymb}
\usepackage{amsthm}
\usepackage{amsfonts}
\usepackage{latexsym}
\usepackage{epsf}

\newcommand{\seqnum}[1]{\href{http://www.research.att.com/cgi-bin/access.cgi/as/~njas/sequences/eisA.cgi?Anum=#1}{\underline{#1}}}

\begin{document}

\begin{center}
\vskip 1cm{\LARGE\bf A characterization of all equilateral
triangles in $\Bbb Z^3$}

\vskip 1cm \large Ray Chandler and  Eugen J. Ionascu\\
\vskip 1cm

 {\small Associate Editor of OEIS,
\href{mailto:RayChandler@alumni.tcu.edu}{\tt
RayChandler@alumni.tcu.edu}\\

and\\

Department of Mathematics\\
Columbus State University\\
Columbus, GA 31907, US \\
{\it Honorific Member of the \\ Romanian Institute of Mathematics
``Simion Stoilow"} \\
\href{mailto:ionascu_eugen@colstate.edu}{\tt ionascu\_eugen@colstate.edu}} \\
\end{center}

\vskip .2in

\begin{abstract}
This paper is a continuation of the work started in \cite{eji1}
and \cite{eji2}. We extend one of the theorems that gave a way to
construct equilateral triangles whose vertices have integer
coordinates to the general situation. An  approximate
extrapolation formula for the sequence $ET(n)$ of all equilateral
triangles with vertices in $\{0,1,2,...,n\}^3$ (A 102698) is given
and the asymptotic behavior of this sequence is analyzed.
\end{abstract}

\newtheorem{theorem}{Theorem}[section]
\newtheorem{proposition}{Proposition}[section]
\newtheorem{corollary}{Corollary}[section]
\newtheorem{lemma}{Lemma}[section]
\newtheorem{definition}{Definition}[section]
\newtheorem{problem}{Problem}[section]
\section{Introduction}

It turns out that equilateral triangles in $\Bbb Z^3$ exist and
there are unexpectedly many. Just to give an example, if we
restrict our attention only to the cube $\{0,1,2,...,1105\}^3$ we
have 2,474,524,936,846,512 of them. In \cite{eji1} it was shown
the first part of the following theorem and the second part about
the converse was only proven under the hypothesis that
$gcd(d,a)=1$ or $gcd(d,b)=1$ or $gcd(d,c)=1$. The main result of
this paper is to show that one can drop this condition.

\begin{theorem}\label{generalpar} Let $a$, $b$, $c$, $d$ be odd positive
integers such that $a^2+b^2+c^2=3d^2$ and $gcd(a,b,c)=1$. Then  the
points $P(u,v,w)$ and $Q(x,y,z)$  whose coordinates given by
\begin{equation}\label{paramone}
\begin{cases}
u=m_um-n_un,\\
v=m_vm-n_vn,\\
w=m_wm-n_wn, \\
\end{cases}
\ \ \ \text{and} \ \ \
\begin{cases}
x=m_xm-n_xn,\\
y=m_ym-n_yn,\\
z=m_zm-n_zn,
\end{cases}\ \
\end{equation}
with
\begin{equation}\label{paramtwo}
\begin{array}{l}
\begin{cases}
m_x=-\frac{1}{2}[db(3r+s)+ac(r-s)]/q,\ \ \ n_x=-(rac+dbs)/q\\
m_y=\frac{1}{2}[da(3r+s)-bc(r-s)]/q,\ \ \ \ \  n_y=(das-bcr)/q\\
m_z=(r-s)/2,\ \ \ \ \ \ \ \ \ \ \ \ \ \ \ \  \ \ \ \ \  \ \ \ \ \ \ n_z=r\\
\end{cases}\\
\ \ \ \text{and}\ \ \\
\begin{cases}
m_u=-(rac+dbs)/q,\ \ \ n_u=-\frac{1}{2}[db(s-3r)+ac(r+s)]/q\\
m_v=(das-rbc)/q,\ \ \ \ \   n_v=\frac{1}{2}[da(s-3r)-bc(r+s)]/q\\
m_w=r,\ \ \ \ \ \ \ \ \ \ \ \ \ \ \ \ \ \ \  n_w=(r+s)/2\\
\end{cases}
\end{array}
\end{equation}
where $q=a^2+b^2$ and $(r,s)$ is a suitable solution of
$2q=s^2+3r^2$ which makes all the numbers in (\ref{paramtwo})
integers, together with the origin (O(0,0,0)) forms an equilateral
triangle in $\Bbb Z^3$ contained in the plane
$${\cal
P}_{a,b,c}:=\{(\alpha,\beta,\gamma)|a\alpha+b\beta+c\gamma=0\}$$
and having sides-lengths equal to $d\sqrt{2(m^2-mn+n^2)}$.

Conversely, there exist a choice of the integers $r$ and $s$ such
that given an arbitrary equilateral triangle contained in the
plane ${\cal P}_{a,b,c}$ with one of the vertices the origin, the
other two vertices are of the form (\ref{paramtwo}) for some
integer values $m$ and $n$.
\end{theorem}

The condition

\begin{equation}\label{degenerate}
min(gcd(d,a),gcd(d,b),gcd(d,c))>1
\end{equation}

defines $a$, $b$, $c$, and $d$ as a {\it degenerate solution} of the
Diophantine equation

\begin{equation}\label{main}
a^2+b^2+c^2=3d^2.
\end{equation}

The first $d$ that admits such degenerate decompositions in $a$,
$b$, and $c$ is $d=1105$. There exactly seven of them:

\[
\begin{array}{c}
(a,b,c)\in \{ (731, 1183, 1315), (475, 1309, 1313), (299, 493,
1825), \\ \\ (1027, 1139, 1145), (187, 415, 1859), (265, 533, 1819),
(493, 1001, 1555)\}
\end{array}
\]

Our study started with the intent of computing the sequence
$ET(n)$ of all equilateral triangles with vertices in
$\{0,1,2,...,n\}$. These values were calculated by the first
author for $n\le 1105$ using a improved version of the code
published in \cite{eji2} and translated into Mathematica (see
A102698 in The On-Line Encyclopedia of Integer Sequences).


One of the parametrizations like in  (\ref{paramtwo}) in the case
$a=731$, $b=1183$ and $c=1315$ is shown below:

\[
\begin{array}{l}
P=(901m-1428n, -1157m+221n, 540m+595n), \\ Q=(-527m-901n,
-936m+1157n, 1135m-540n), m,n\in \Bbb Z .
\end{array}
\]

For other values of $r$ and $s$ we get formally different
parametrizations but they are all equivalent in the sense that
they can be obtained from the above by changing the variables in
the expected ways:

\begin{equation}\label{hexagonaltransform}
\begin{array}{c}
(m,n)\to (-m,-n), \ (m,n)\to (m-n,m), \ (m,n)\to (n-m,n),\\ \\
(m,n)\to (m-n,-n) \ and \ (m,n)\to (n-m,-m). \end{array}
\end{equation}

 These changes of variables leave invariant the quadratic form involved in the
side-length of the triangle the various points given by
(\ref{paramone}) define the vertices of a tessellation of the
plane ${\cal P}_{a,b,c}$ with equilateral triangles.

The above example is saying and it is supporting evidence that the
parametrization (\ref{paramtwo}) is essentially  the only one and
the Theorem~\ref{generalpar} can be extended to cover the
degenerate solutions case also.

\section{The arguments of the Theorem~\ref{generalpar}. }

\begin{proof} The first part of the Theorem~\ref{generalpar} follows from
\cite{eji1}. For the second part we are going to reconstitute some
of the details started in the proof of the particular case:
$gcd(d,c)=1$.

Let us start with a triangle in ${\cal P}_{a,b,c}$ say $\triangle
OPQ$ with $P(u_0,v_0,w_0)$ and $Q(x_0,y_0,z_0)$.

By Theorem~4 in \cite{eji1}we have

\begin{equation}\label{eqfortp}
\begin{cases}
\displaystyle x_0=\frac{u_0}{2}\pm \frac{cv_0-bw_0}{2d}\\
\displaystyle y_0=\frac{v_0}{2}\pm \frac{aw_0 -cu_0}{2d}  \\
\displaystyle z_0 =\frac{w_0}{2}\pm \frac{bu_0-av_0}{2d},
\end{cases}
\end{equation}

\noindent for some choice of the signs.  This means that $d$ must
divide $cv_0-bw_0$, $aw_0-cu_0$ and $bu_0-aw_0$.

So, we need to look at the the following system of linear
equations in $m$ and $n$:

\begin{equation}\label{fsystem}
\begin{cases}
u_0=m_um-n_un,\\
v_0=m_vm-n_vn,\\
w_0=m_wm-n_wn. \\
\end{cases}
\end{equation}

By the Kronecker-Capelli theorem this linear system of equations has
a solution if and only the rank of the main matrix is the same as
the rank of the extended matrix. Since $m_vn_w-m_wn_v=ad$,
$m_wn_u-m_un_w=bd$ and $m_un_v-m_vn_u=cd$ (one can check these
calculations based on the definitions in (\ref{paramtwo})) then the
rank of the main matrix is two and the rank of the extended matrix
is also two because its determinant is $u_0ad+v_0bd+w_0dc=0$.
\noindent  This implies  that (\ref{fsystem}) has a unique real (in
fact rational) solution in $m$ and $n$.

We want to show that this solution is in fact an integer solution.
Solving for $n$ from each pair of equations in (\ref{fsystem}) we
get

\begin{equation}\label{n1}
n=\frac{v_0m_w-w_0m_v}{ad}=\frac{w_0m_u-u_0m_w}{bd}=\frac{u_0m_v-v_0m_u}{cd}.
\end{equation}

Because $gcd(a,b,c)=1$, there exist integers $a'$, $b'$, $c'$ such
that $aa'+bb'+cc'=1$. Then one can see that

\begin{equation}\label{n2}
n= \frac{a'(v_0m_w-w_0m_v)+b'(w_0m_u-u_0m_w)+c'(u_0m_v-v_0m_u)}{d}
\end{equation}

Next, from (\ref{n2}) we observe that in order for $n$ to be an
integer it is enough to prove that $d$ divides $v_0m_w-w_0m_v$,
$w_0m_u-u_0m_w$ and $u_0m_v-v_0m_u$. Hence we calculate for example
$v_0m_w-w_0m_v$ in more detail:
 $$\begin{array}{l}
 \displaystyle v_0m_w-w_0m_v=v_0r-\frac{das-rbc}{q}w_0=\frac{v_0qr-(das-rbc)w_0}{q}=\\ \\
 \displaystyle \frac{-dasw_0+v_0(3d^2-c^2)r+rbcw_0}{q}=\frac{c(bw_0-v_0c)r+3rv_0d^2-dasw_0}{q}.
\end{array}$$
From (\ref{eqfortp}) we see that $bw_0-v_0c=\pm d(2x_0-u_0)$. Hence,

\begin{equation}\label{finalstep}
\displaystyle v_0m_w-w_0m_v=\frac{d[\pm
c(2x_0-u_0)r+3rv_0d-asw_0]}{q}.
\end{equation}

Assuming that $gcd(d,c)=\zeta$ we can write $d=\zeta d_1$ and
$c=\zeta c_1$ with $gcd(d_1,c_1)=1$. Also we see that $\zeta ^2$
must divide $q=3d^2-c^2$ so let us write $q=\zeta^2 q_1$. If $p$
is a prime dividing $\zeta$, it must be an odd prime and if it is
of the form $4k+3$ it must divide $a$ and $b$ which is
contradicting the assumption that $gcd(a,b,c)=1$. Therefore it
must be a prime of the form $4k+1$. Hence $q_1$ is still a sum of
two squares.

In the proof of Theorem~13 in \cite{eji1} one can choose $r$ and
$s$ with the extra condition that $r$ and $s$ are divisible by
$\zeta$. Indeed, Lemma~14 in \cite{eji1} is applied to
$(ac)^2+3(db)^2=\zeta^2[(ac_1)^2+3(d_1b)^2]$ and to $q=\zeta^2
q_1$ but instead one can apply it to $(ac_1)^2+3(d_1b)^2$ and to
$q_1$ giving, let us, say $r_1$ and $s_1$. Then  we put $r=\zeta
r_1$ and $s=\zeta s_1$ and then all the arguments there go as
stated. From (\ref{finalstep}) we see that

\begin{equation}\label{finalstep2}
\begin{array}{c}
\displaystyle v_0m_w-w_0m_v=\\ \\  \displaystyle \frac{d[\pm
c(2x_0-u_0)r+3rv_0d-asw_0]}{q}=  \frac{\zeta d_1 [\pm \zeta
c_1(2x_0-u_0)\zeta
r_1+3\zeta r_1 v_0\zeta d_1-a\zeta s_1 w_0]}{\zeta ^2 q_1}=\\ \\
\displaystyle \frac{\zeta^2 d_1[ \pm c_1(2x_0-u_0)r_1\zeta+3r_1
v_0d_1\zeta -a s_1 w_0]}{\zeta ^2 q_1}= \frac{ d_1 \xi}{ q_1},
\end{array}
\end{equation}

\noindent where $\xi= \pm c_1(2x_0-u_0)r_1\zeta+3r_1v_0d_1\zeta
-as_1 w_0$. This implies that $d_1$ must divide $v_0m_w-w_0m_v$.
In a similar way we can show that $d_1$ divides $w_0m_u-u_0m_w$
and $u_0m_v-v_0m_u$ and so from (\ref{n2}) we see that $n$ is a
rational with denominator having $\zeta$ as a multiple. Similar
arguments will give us that $m$ is of the same form.

The triangle having the coordinates as in (\ref{paramone}) with
these $m$ and $n$ will give the length

$$l^2=2d^2(m^2-mn+n^2).$$

But this whole construction can be repeated for $a$ or $b$ instead
of $c$ and we obtain that

$$l^2=2d^2(m_1^2-m_1n_1+n_1^2),$$

\noindent for some rational numbers $m_1$, $n_1$ with denominator
having multiple  $\eta=gcd(d,b)$ for example.  Since
$gcd(\zeta,\eta)=1$ we see that $m^2-mn+n^2=m_1^2-m_1n_1+n_1^2$ must
be an integer. Therefore

\begin{equation}\label{important}
l^2=2d^2(\alpha^2-\alpha\beta+\beta^2), \alpha,\beta\in \Bbb Z.
\end{equation}

On the other hand if $u_0^2+v_0^2+w_0^2=l^2$ and $\zeta$ divides $d$
we can see that

\begin{equation}\label{key}
u_0^2+v_0^2+w_0^2\equiv 0\ (mod \ \zeta^2) .
\end{equation}

We also know that $au_0+bv_0+cw_0=0$ and hence
$a^2u_0^2=b^2v_0^2+2bcv_0w_0+c^2w_0^2$ . This implies
$a^2u_0^2\equiv b^2v_0^2+2bcv_0w_0$ (mod $\zeta^2$). But
$a^2+b^2=3d^2-c^2\equiv 0$ (mod $\zeta^2$) too and then
$a^2(u_0^2+v_0^2)\equiv 2bcv_0w_0$ (mod $\zeta^2$) which combined
with (\ref{key}) gives

\begin{equation}\label{last}
a^2w_0^2+2bcv_0w_0\equiv 0\ (mod \ \zeta^2).
\end{equation}

Because we must have $gcd(a,\zeta)=1$, (\ref{last}) implies that
$\zeta$ divides $w_0$.  Indeed, if $p$ is a prime that has
exponent one in the decomposition of $\zeta$ then (\ref{last})
gives in particular $w_0^2\equiv 0$ (mod $p$) and so $p$ must
divide $w_0$. If the exponent of $p$ in $\zeta$ is two, then
(\ref{last}) in particular implies that $w_0$ is divisible by $p$
but then $aw_0^2\equiv 0$ (mod $p^3$) which implies $p^2$ divides
$w_0$. Inductively if the exponent of $p$ in $\zeta$ is $k$ then
this must be true for $w_0$ too. Hence we must have $\zeta$ a
divisor of $w_0$ so $w_0=\zeta w_0'$.

Now we can go back to (\ref{finalstep2}) and observe that we can
rewrite it as

\begin{equation}\label{back}
v_0r-w_0m_v=\frac{\zeta d_1\xi'}{q_1},
\end{equation}
\noindent where $\xi'= \pm c_1(2x_0-u_0)r_1+3r_1v_0d_1 -as_1
w_0'$. Now we observe that the left hand side of (\ref{back}) is a
multiple of $\zeta$ since $r$ and $w_0$ are. After simplification
with $\zeta$ this will show that $q_1$ must divide in fact $\xi'$
and so $d$ divides $v_0m_w-w_0m_v$. Similar arguments can be used
to deal with the other cases. Hence $n$ must be an integer and so
should be $m$. Changing the variables as in
(\ref{hexagonaltransform}),  one of the corresponding triangles
given by (\ref{paramone}) is going to match with the triangle
$OPQ$.\end{proof}

{\bf Remark:} One can see that the condition on $r$ and $s$ to be
divisible by $\zeta$ is implied by asking only that the numbers in
(\ref{paramtwo}) be integers. Indeed, given such a choice of $r$
and $s$, they will define by (\ref{paramone}), in which $m=1$ and
$n=0$, an equilateral triangle with integer coordinates. According
to the proof of Theorem~\ref{generalpar}, $w=r$ and $z=(r-s)/2$
must be divisible by $\zeta$. This implies that $r$ and $s$ must
be multiples of $\zeta$. Therefore any parametrization as in the
Theorem~\ref{generalpar} is unique up to the transformations
(\ref{hexagonaltransform}).

\section{Behavior of the sequence $ET(n)$ }
\begin{figure}
\begin{center}
\epsfig{file=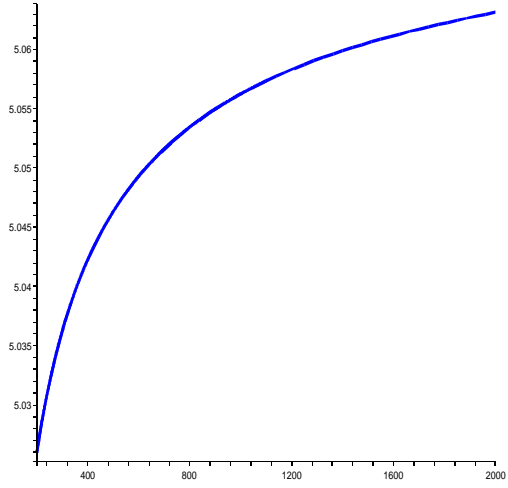,height=3in,width=3in} \caption{The graph of
$g$ extrapolating $f$ over the interval $[1,2000]$} \label{figure1}
\end{center}
\end{figure}

The calculations of the $ET(n)$ for all $n\le 1105$ gave us enough
data to be able to extrapolate the graph of $n\overset{f}{\to}
\frac{\ln (ET(n))}{\ln (n+1)}$ as shown in Figure~\ref{figure1}. The
function we used to extrapolate is of the form
$g(x)=a+\frac{b}{\sqrt{x}+c}$ having clearly $a$ as limit at
infinity. Then we made it agree with $f$ on three points. That gave
us $a:=5.079282921$, $b:=-0.7091588389$, and $c:=-0.8403164433$.
Numerically then we discovered that the average of $|f(k)-g(k)|$
over all values of $k=1,...,1105$ is approximately $0.002638971108$.

One conjecture that we would like to make here is that $f(n)$ is a
strictly increasing sequence and then as result it is convergent to
a constant $C\approx 5.08$.

\begin{figure}
\begin{center}
\epsfig{file=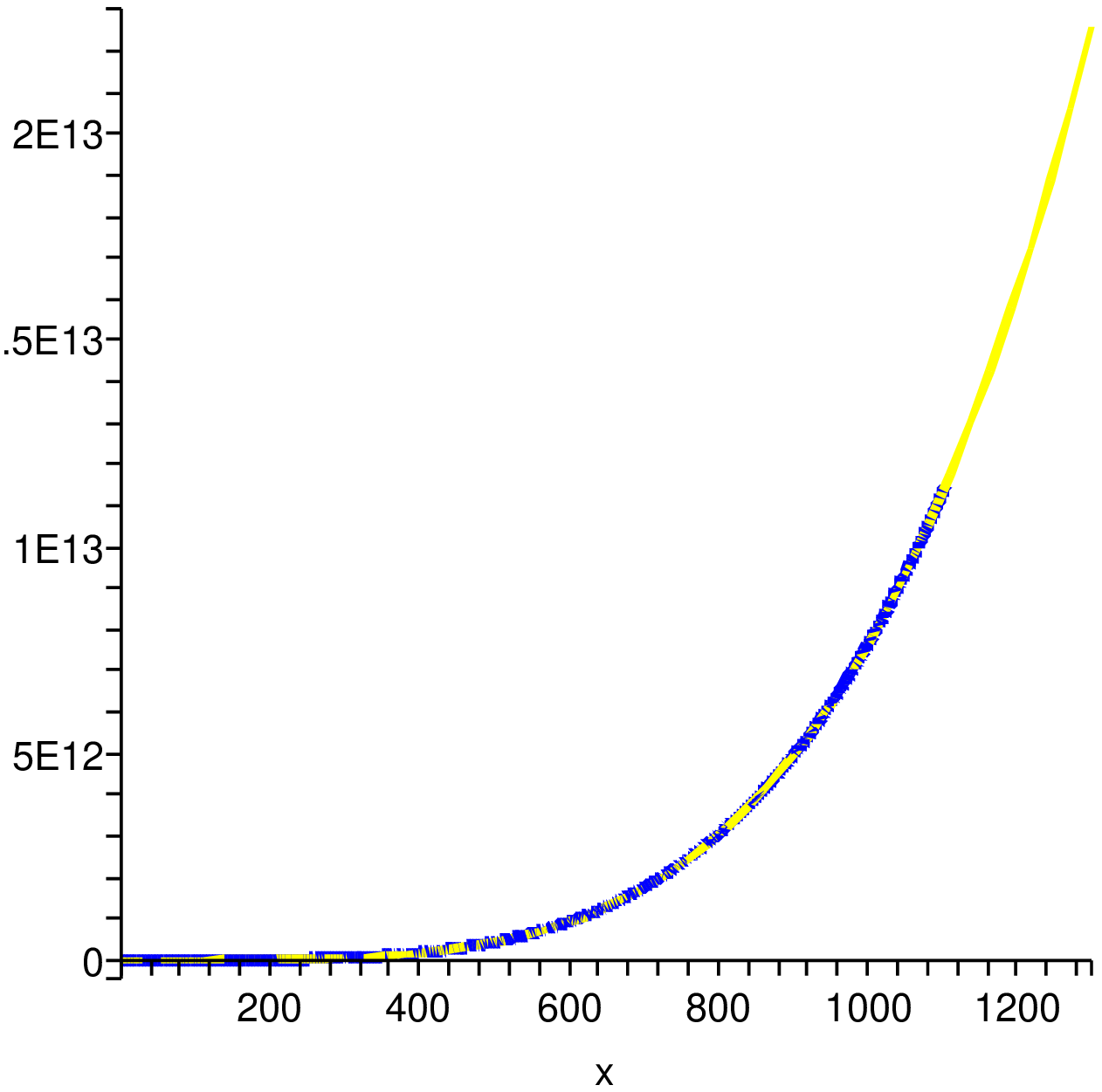,height=3in,width=3in}
  \caption{The graph of $n\to ET(n+1)-ET(n)$}
\label{derivative}
\end{center}
\end{figure}

\begin{figure}
\begin{center}
\epsfig{file=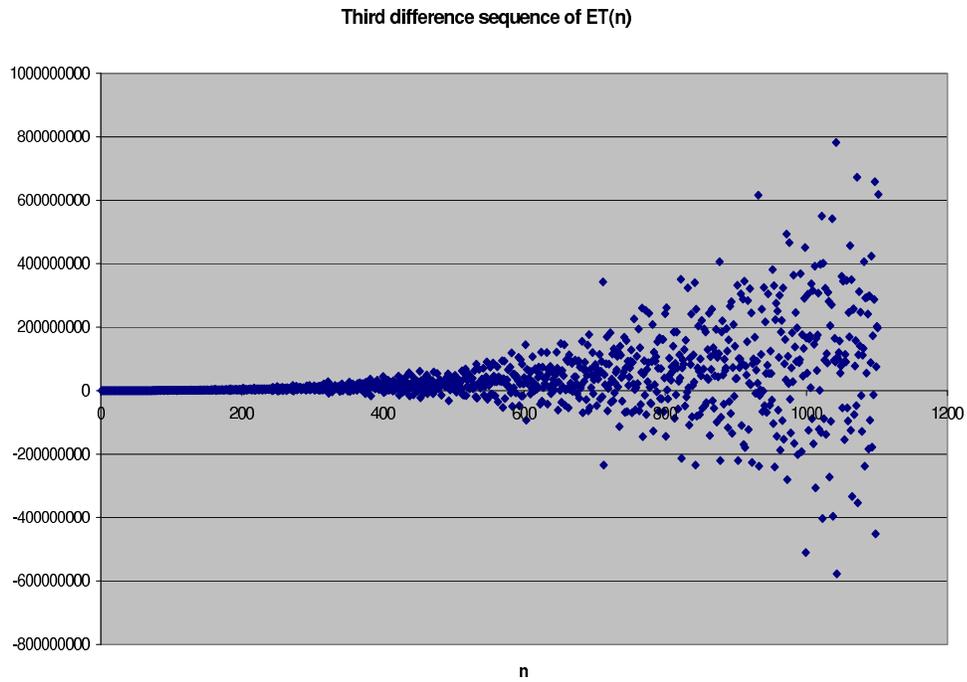,height=3.5in,width=5in} \caption{$\Delta ^3
ET(n)$}

\label{3rdder}
\end{center}
\end{figure}

The graph of the ``derivative" of $ET(n)$ (Figure~\ref{derivative})
is almost like the graph of $h(x)=C(x+1)^k$ where $k:=4.151431798$
and $C:=2.660972140$. The third difference of $ET(n)$ as represented
in Figure~\ref{3rdder} seems to bring a chaotic flavor to this
sequence and it is saying in a certain sense that no simple formula
for $ET(n)$ can exist.

\eject

\bigskip
\hrule
\bigskip

\noindent 2000 {\it Mathematics Subject Classification}: Primary
11A67;
Secondary 11D09, 11D04, 11R99, 11B99, 51N20.\\

\noindent {\it Keywords}: Diophantine equations, equilateral
triangles, integers, parametrization, characterization.

\bigskip
\hrule
\bigskip

\noindent (Concerned with sequence \seqnum{A102698}.)

\bigskip
\hrule
\bigskip




\end{document}